\documentclass [12pt,a4paper,reqno]{amsart}
\input amssymb.sty
\usepackage{amsmath,amsthm}
\usepackage{amsfonts}
\usepackage{amssymb}
\allowdisplaybreaks

\newcommand{\be}{\begin{equation}}
\newcommand{\ef}{\end{equation}}
\chardef\bslash=`\\ 

\newtheorem{thm}{Theorem}[section]
\newtheorem*{thm*}{Theorem}
\newtheorem{cor}[thm]{Corollary}
\newtheorem{lem}[thm]{Lemma}
\newtheorem{prop}[thm]{Proposition}

\theoremstyle{definition}

\newtheorem*{remark*}{Remarks}
\newtheorem*{defn*}{Definition}

\theoremstyle{remark}

\numberwithin{equation}{section}

\newcommand{\wt}{\widetilde}
\newcommand{\wh}{\widehat}
\newcommand{\fc}{\frac}
\newcommand{\bk}{\bigskip}
\newcommand{\iy}{\infty}
 \renewcommand{\sectionmark}[1]{}

\newcommand{\Be}{Beltrami}
\newcommand{\Gr} {Grunsky}
\newcommand{\hol} {holomorphic}
\newcommand{\qc} {quasiconformal}
\newcommand{\sh} {subharmonic}
\newcommand{\psh} {plurisubharmonic}

\newcommand{\Te} {Teichm\"{u}ller}
\newcommand{\Ko} {Kobayashi}
\newcommand{\Ca} {Carath\'{e}odory}
\newcommand{\uTs} {universal Teichm\"{u}ller space}
\newcommand{\const}{\operatorname{const}}

 \usepackage{amsfonts}
\newcommand{\field}[1]{\mathbb{#1}}

\newcommand{\dl}{\delta}
\newcommand{\D}{\Delta}
\newcommand{\om}{\omega}
\newcommand{\z}{\zeta}
\newcommand{\ov}{\overline}
\newcommand{\vp}{\varphi}
\newcommand{\hC}{\wh{\field{C}}}
\newcommand{\C}{\field{C}}

\newcommand{\B}{\mathbf{B}}
\newcommand{\T}{\mathbf{T}}
\newcommand{\Belt}{\mathbf{Belt}}

\newcommand{\dist}{\operatorname{dist}}

\newcommand{\Om} {\Omega}
\newcommand{\vk} {\varkappa}
\newcommand{\x} {\mathbf x}
\renewcommand{\a} {\alpha}

\newcommand{\ld}{\lambda}

\newcommand{\kp}{\kappa}

\begin{document}

\title{Strengthened Grunsky and Milin inequalities}
\author{Samuel L. Krushkal}

\begin{abstract}
The method of \Gr \ inequalities has many applications and has been 
extended in many directions, even to bordered Riemann surfaces.
However, unlike the case of functions univalent in the disk, a \qc
\ variant of this theory has not been developed so far. In this
paper, we essentially improve the basic facts concerning the
classical \Gr \ inequalities for univalent functions on the disk
and extend these results to arbitrary quasiconformal disks. 
Several applications are given. 
\end{abstract}

\date{\today\hskip4mm({grumil.tex})}

\maketitle

\bigskip

{\small {\textbf {2010 Mathematics Subject Classification:}
Primary: 30C55, 30C62, 31A35; Secondary 30F60, 32F45}

\medskip


\bigskip

\markboth{Samuel L. Krushkal}{Strengthened Grunsky and Milin
inequalities} \pagestyle{headings}

\bk 
\section{The Grunsky and Grunsky-Milin coefficients}

\subsection{The Grunsky operator}
In 1939, H. \Gr \ discovered the necessary and sufficient conditions
for univalence of a \hol \ function in a finitely connected domain
on the extended complex plane $\hC = \C \cup \{\iy\}$ in terms of
an infinite system of the coefficient inequalities. In particular,
his theorem for the canonical disk $\D^* = \{z \in \hC: \ |z| >
1\}$ yields that a \hol \ function $f(z) = z + \const + O(z^{-1})$
in a neighborhood of $z = \iy$ can be extended to a univalent \hol
\ function on the $\D^*$ if and only if its \Gr \ coefficients
$\a_{mn}$ satisfy \be\label{1.1} \Big\vert \sum\limits_{m,n=1}^\iy
\ \sqrt{m n} \ \a_{mn} x_m x_n \Big\vert \le 1,
\end{equation}
where $\a_{mn}$ are defined by \be\label{1.2} \log \fc{f(z) -
f(\z)}{z - \z} = - \sum\limits_{m,n=1}^\iy \a_{mn} z^{-m} \z^{-n},
\quad (z, \z) \in (\D^*)^2,
\end{equation} 
the sequence $\mathbf x = (x_n)$ runs over the unit sphere 
$S(l^2)$ of the Hilbert space $l^2$ with 
norm $\|\x\|^2 = \sum\limits_1^\iy |x_n|^2$,
and the principal branch of the logarithmic function is chosen
(cf. \cite{Gr}). The quantity
\be\label{1.3} 
\vk(f) = \sup \Big\{\Big\vert \sum\limits_{m,n=1}^\iy \ 
\sqrt{mn} \ \a_{mn} x_m x_n \Big\vert: \   
\mathbf x = (x_n) \in S(l^2) \Big\} \le 1
\end{equation}
is called the \textbf{Grunsky norm} of $f$.

For the functions with $k$-\qc \ extensions ($k < 1$), we have 
instead of (1.3) a stronger bound 
\be\label{1.4} 
\Big\vert \sum\limits_{m,n=1}^\iy \ \sqrt{mn} \ \a_{mn} x_m x_n 
\Big\vert \le k \quad \text{for any} \ \ \x = (x_n) \in S(l^2), 
\end{equation}
established first in \cite{Ku1} (see also \cite{Kr7}). 
Then $\vk(f) \le k(f)$, where $k(f)$ denotes the \textbf{\Te \
norm} of $f$ which is equal to the infimum of dilatations 
$k(w^\mu)= \|\mu\|_\iy$ of \qc \ extensions of $f$ to $\hC$. 
Here $w^\mu$ denotes a homeomorphic solution to the \Be \ equation
$\partial_{\ov z} w = \mu \partial_z w$ on $\C$ extending $f$; 
accordingly, 
$\mu$ is called the \textbf{\Be \ coefficient} (or complex dilatation) 
of $w$.

Note that the \Gr \ (matrix) operator 
$\mathcal G(f) = (\sqrt{m n} \ \a_{mn}(f))_{m,n=1}^\iy$ 
acts as a linear operator $l^2 \to l^2$
contracting the norms of elements $\x \in l^2$; the norm of this
operator equals $\vk(f)$.

For most functions $f$, we have the strong inequality $\vk(f)
< k(f)$, while the functions with the equal norms play a crucial
role in many applications.

\subsection{Generalization}
The method of \Gr \ inequalities was generalized in several
directions, even to bordered Riemann surfaces $X$ with a finite
number of boundary components(cf. \cite{Gr}, \cite{Le}, \cite{Mi},
\cite{Po}, \cite{SS}). In the general case, the generating
function (1.2) must be replaced by a bilinear differential
\be\label{1.5} 
- \log \fc{f(z) - f(\z)}{z - \z} - R_X(z, \z) =
\sum\limits_{m, n = 1}^\iy \beta_{m n} \ \vp_m(z) \vp_n(\z): \ X
\times X \to \C,
\end{equation}
where the surface kernel $R_X(z, \z)$ relates to the conformal map
$j_\theta(z, \z)$ of $X$ onto the sphere $\hC$ slit along arcs of
logarithmic spirals inclined at the angle $\theta \in [0, \pi)$ to
a ray issuing from the origin so that $j_\theta(\z, \z) = 0$ and
$$
j_\theta(z) = (z - z_\theta)^{-1} + \const + O(1/(z - z_\theta)) \quad
\text{as} \ \ z \to z_\theta = j_\theta^{-1}(\iy)
$$
(in fact, only the maps $j_0$ and $j_{\pi/2}$ are applied). Here
$\{\vp_n\}_1^\iy$ is a canonical system of \hol \ functions on $X$
such that (in a local parameter)
$$
\vp_n(z) = \fc{a_{n,n}}{z^n} + \fc{a_{n+1,n}}{z^{n+1}} + \dots
\quad \text{with} \ \ a_{n,n} > 0, \quad n = 1, 2, \dots ,
$$
and the derivatives (linear \hol \ differentials) $\vp_n^\prime$
form a complete orthonormal system in $H^2(X)$.

We shall deal only with simply connected domains $X = D^* \ni \iy$
with \qc \ boundaries (quasidisks). For any such domain, the
kernel $R_D$ vanishes identically on $D^* \times D^*$, and the
expansion (1.5) assumes the form
 \be\label{1.6}
 - \log \fc{f(z) - f(\z)}{z - \z} 
= \sum\limits_{m, n = 1}^\iy 
\fc{\beta_{m n}}{\sqrt{m n} \ \chi(z)^m \ \chi(\z)^n},
\end{equation}
where $\chi$ denotes a conformal map of $D^*$ onto the disk $\D^*$
so that $\chi(\iy) = \iy, \ \chi^\prime(\iy) > 0$. 

Each coefficient $\beta_{m n}(f)$ in (1.6) is represented as a
polynomial of a finite number of the initial coefficients $b_1,
b_2, \dots, b_s$ of $f$; hence it depends holomorphically on \Be \
coefficients of \qc \ extensions of $f$ as well as on the
Schwarzian derivatives 
\be\label{1.7} 
S_f(z) = \Bigl(\fc{f^{\prime\prime}(z)}{f^\prime(z)}\Bigr)^\prime -
\fc{1}{2} \Bigl(\fc{f^{\prime\prime}(z)}{f^\prime(z)}\Bigr)^2,
\quad z \in D^*. 
\end{equation}
These derivatives range over a bounded domain in the complex 
Banach space $\B(D^*)$ of hyperbolically bounded \hol \ functions 
$\vp \in \D^*$ with norm 
$$
\|\vp\|_\B = \sup_{D^*} \ld_{D^*}^{-2}(z) |\vp(z)|, 
$$
where $\ld_{D^*}(z)|dz|$ denotes the hyperbolic metric of $D^*$ 
of Gaussian curvature $- 4$. 
This domain models the \textbf{\uTs} \ $\T$ with 
the base point $\chi^\prime(\iy) D^*$ (in \hol \ Bers' embedding 
of $\T$). 

A theorem of Milin extending the Grunsky univalence criterion for
the disk $\D^*$ to multiply connected domains $D^*$ states that a
\hol \ function $f(z) = z + \const + O(z^{-1})$ in a neighborhood
of $z = \iy$ can be continued to a univalent function in the whole
domain $D^*$ if and only if the coefficients $\a_{mn}$ in (1.6) 
satisfy, 
similar to the classical case of the disk $D^*$, the inequality
\be\label{1.8} 
\Big\vert \sum\limits_{m,n=1}^\iy \ \beta_{m n} \ x_m x_n 
\Big\vert \le 1
\end{equation}
for any point $\x = (x_n) \in S(l^2)$. We call the quantity
\be\label{1.9} 
\vk_{D^*}(f) = \sup \Big\{ \Big\vert
\sum\limits_{m,n = 1}^{\iy} \ \beta_{mn} \ x_m x_n \Big\vert : \
{\mathbf x} = (x_n) \in S(l^2)\Big\},
\end{equation}
the \textbf{generalized \Gr \ norm} of $f$.

Note that in the case 
$D^* = \D^*$, \ $\beta_{mn} = \sqrt{m n} \ \a_{mn}$; 
for this disk, we shall use the notations $\Sigma$ and
$\vk(f)$.

By (1.8), $\vk_{D^*}(f) \le 1$ for any $f$ from the class
$\Sigma(D^*)$ of univalent functions in $D^*$ with hydrodynamical
normalization
$$
f(z) = z + b_0 + b_1 z^{-1} + \dots \quad \text{near} \ \ z = \iy.
$$
However, unlike the case of functions univalent in the disk, a \qc
\ variant of this theory has not been developed so far. 

\subsection{}
The technique of the \Gr \ inequalities is a powerful tool in 
geometric complex analysis having 
fundamental applications in the \Te \ space theory and other
fields and concerns mainly the classical case of univalent
functions on the disk $\D^*$ with hydrodynamical normalization,
which has been investigated by many authors from different points
of view. 

In this paper, we create the \qc \ theory of generic \Gr \ coefficients 
and essentially improve the basic facts and estimates 
concerning the classical \Gr \ inequalities. These results are 
extended to univalent functions on arbitrary \qc \ disks.

\section{Main results}

\subsection{}
First recall the fundamental property of extremal Beltrami
coefficients which plays a crucial role in applications of
univalent functions with \qc \ extensions. Consider the unit ball
of \Be \ coefficients
$$
\Belt(D)_1 = \{\mu \in L_\iy(C): \ \mu(z)|D^* = 0, \ \ \|\mu\|_\iy
< 1\}
$$
and their pairing with $\psi \in L_1(D)$ by
$$
\langle \mu, \psi\rangle_D = \iint\limits_D  \mu(z) \psi(z) dx dy
\quad (z = x + i y).
$$
The following two sets of \hol \ functions $\psi$ (equivalently,
of \hol \ quadratic differentials $\psi dz^2$)
$$
\begin{aligned}
A_1(D) &= \{\psi \in L_1(D): \ \psi \ \ \text{\hol \ in} \ \ D\},  \\
A_1^2(D) &= \{\psi = \om^2 \in A_1(D): \ \om \ \ \text{\hol \ in}
\ \ D\}
\end{aligned}
$$
are intrinsically connected with the extremal \Be \ coefficients
(hence, with the \Te \ norm) and \Gr \ inequalities.
The well-known criterion for extremality (the
Hamilton-Krushkal-Reich-Strebel theorem) implies that a \Be \
coefficient $\mu_0 \in \Belt(D)_1$ is extremal if an only if
\be\label{2.1} \|\mu_0\|_\iy = \sup_{\|\psi\|_{A_1(D)}=1} |\langle
\mu_0, \psi \rangle_D |.
\end{equation}
The same condition is necessary and sufficient for the
infinitesimal extremality of $\mu_0$ (i.e., at the origin of $\T$
in the direction $t \phi_\T(\mu_0)$, where $\phi_\T$ is the
defining (factorizing) \hol \ projection $\Belt(D)_1 \to \T$ );
see, e.g., \cite{EKK}, \cite{GL}. In contrast, the \Gr \ norm
relates to the functions from $A_1^2(D)$, i.e. to abelian
differentials.

For an element $\mu \in \Belt(D)_1$ we define
$$
\mu^*(z) =\mu(z)/\|\mu\|_\iy,
$$
so that $\|\mu^*\|_\iy = 1$, and associate with the corresponding
map $f^\mu$ the quantity 
\be\label{2.2} 
\a_D(f^\mu) = \sup \
\Big\{ \Big\vert \iint\limits_D \mu^*(z) \vp(z) dx dy \Big\vert: \
\vp \in A_1^2(D), \ \|\vp\|_{A_1} = 1 \Big\} \le 1.
\end{equation}
For the disk $D = \D$, we shall use the notation $\a(f^\mu)$.

\subsection{Strengthened bounds for \Gr \ norm} 
Now we can formulate our results.  The following theorem
essentially improves the basic estimate (1.4).

\begin{thm} For any quasidisk $D^*$, the generalized \Gr \ norm
$\vk_{D^*}(f)$ of every function $f \in \Sigma^0(D^*)$ is
estimated by its \Te \ norm $k = k(f)$ by 
\be\label{2.3}
\vk_{D^*}(f) \le k \fc{k + \a_D(f)}{1 + \a_D(f) k},
\end{equation}
and $\vk_{D^*}(f) < k$ unless $\a_D(f) = 1$. 
The last equality occurs if and only if $\vk_{D^*}(f) = k(f)$. 
\end{thm}

\begin{thm} The equality $\vk_{D^*}f = k(f)$ holds if and
only if the function $f$ is the restriction to $\ov{D^*}$ of 
a \qc \ self-map $w^{\mu_0}$ of $\hC$ with \Be \ coefficient 
$\mu_0$ satisfying the condition
\be\label{2.4}
\sup |\langle \mu_0, \vp\rangle_D| = \|\mu_0\|_\iy,
\end{equation}
where the supremum is taken over \hol \ functions $\vp \in A_1^2(D)$
with $\|\vp\|_{A_1(D)} = 1$.

If, in addition, the equivalence class of $f$ (the collection of
maps equal $f$ on $\partial D^*$) is a Strebel point, then $\mu_0$ is
necessarily of the form 
\be\label{2.5}
\mu_0(z) = \|\mu_0\|_\iy |\psi_0(z)|/\psi_0(z) \quad
\text{with} \ \psi_0 \in A_1^2(D).
\end{equation}
\end{thm} 

The condition (2.4) has a geometric nature based on the properties 
of the invariant \Ca \ and \Ko \ distances of the \uTs \ $\T$. 

The assertion of Theorem 2.2 was earlier established in \cite{Kr2} 
only for the functions univalent in the canonical disk $\D^*$, 
i.e., for $f \in \Sigma$. This special result answered a question 
posed by several mathematicians and has many applications.  

Shiga and Tanigawa gave an essential extension of this phenomena
to \Te \ spaces of elementary groups (see [ShT]). 
In particular, it holds for covers of conformal maps of the punctured 
disk $\{1 < |z| < \iy\}$. 

For $f \in \Sigma$, mapping the unit circle onto an 
analytic curve, the equality (2.5) was obtained by a different 
method in [Ku2].

\subsection{Two corollaries} 
Both Theorems 2.1 and 2.2 have many interesting consequences. 
In this paper we present the consequences of Theorem 2.1. 
We start with corollaries concerning the maps with 
small dilatations. 
   
From (2.3), for all $f \in \Sigma^0(D^*)$ with small dilatation $k(f)$, 
$$ 
\vk_{D*}(f) \le \a_D(f) k + O(k^2),  
$$  
where the bound for the remainder is uniform when $k \le k_0$ 
and $k_0 < 1$ is fixed.   
On the other hand, as was established in \cite{Kr4}, if a 
function $f \in \Sigma^0(D^*)$ admits \qc \ extension $w^\mu$ 
of \Te \ type, i.e. with $\mu = k |\psi|/\psi, \ \psi \in A_1(D)$, then 
its \Gr \ norm is estimated from below by 
\be\label{2.6} \vk_{D*}(f) \ge \a_D(f) k(f),
\end{equation}
with $\a_{D^*}(f)$ given by (2.2). 
Hence, the inequalities (2.3) and (2.6) imply 

\begin{cor} The generalized \Gr \ norm of any $f \in \Sigma^0(D^*)$ 
with \Te \ \qc \ extension 
satisfies the asymptotic equality 
\be\label{2.7} \vk_{D^*}(f) = \a_D(f) k + O(k^2), 
\quad k = k(f) \to 0.
\end{equation}
\end{cor}

\bk In the case of the canonical disk $\D^*$, one obtains from the last
equality a quantitative relation between the \Gr \ norm and the
Schwarzian derivative of $f$. Namely, using the Ahlfors-Weill \qc \ 
extension of univalent functions and letting 
\be\label{2.8}
\nu_\vp(z) = \frac{1}{2}(1 - |z|^2)^2 \vp (1/{\bar z}) 1/{\bar
z}^4, \ \ \vp \in \B = \B(\D^*), 
\end{equation}
one derives

\begin{cor} For $f \in \Sigma(\D^*)$ with sufficiently small norm 
$\|S_f\|_{\B}$ of its Schwarzian, 
\be\label{2.9} 
\vk(f) = \sup \{|\langle \nu_{S_f}, \psi \rangle_{\D}|: 
\ \psi \in A_1^2, \ \|\psi \|_{A_1(\D)} = 1\} + O(\|S_f\|_{\B}^2),
\end{equation}
where the ratio $O(\|S_f\|_{\B}^2)/\|S_f\|_{\B}^2$ remains bounded
as $\|S_f\|_{\B} \to 0$.
\end{cor}

The \Be \ coefficients of the form (2.8) are called harmonic,  
in view of their connection with the deformation theory for 
conformal structures.  

\bk
\subsection{Continuity} 
It is well known that the classical
\Gr \ norm $\vk(f)$ regarded as a curve functional is lower
semicontinuous in the weak topology on the space $\Sigma^0$ (i.e.,
with respect to locally uniform convergence of sequences
$\{f_n\} \subset \Sigma^0$ on the disk $\D^*$) and continuous with
respect to convergence of $f_n$ in \Te \ metric (see \cite{Sc},
\cite{Sh}). The arguments exploited in the proofs essentially use
the univalence on the canonical disk $\D^*$. The continuity of
$\vk(f)$ plays a crucial role in some applications of the \Gr \
inequalities technique to \Te \ spaces.

We consider univalent functions on generic quasidisks $D^*$ and
show that in either case the \Gr \ norm is lower semicontinuous 
in the 
weak topology on $\Sigma^0(D^*)$  and locally Lipschitz
continuious with respect to \Te \ metric.

\begin{thm} (i) \ If a sequence $\{f_n\} \subset \Sigma^0(D^*)$ is
convergent locally uniformly on $D^*$ to $f_0$, then
\be\label{2.10} \vk_{D^*}(f_0) \le \liminf\limits_{n\to \iy}
\vk_{D^*} (f_n).
\end{equation}

(ii) \ The functional $\vk_{D^*}(\vp)$ regarded as a function of
points $\vp = S_f$ from the \uTs \ $\T$ (with base point $D^*$) is
locally Lipschitz continuous and logarithmically \psh \ on $\T$.
\end{thm}

\medskip
This key theorem is essential in the proof of other theorems.

\subsection{Generalization of Moser's conjecture} 
In 1985, J. Moser conjectured that {\em the set of functions $f
\in \Sigma^0$ with $\vk(f) = k(f)$ is rather sparse in $\Sigma^0$
so that any function $f \in \Sigma^0$ is approximated by functions
$f_n$ satisfying $\vk(f_n) < k(f_n)$ uniformly on compact sets in
$\D^*$}. This conjecture was proved in \cite{KK1} and in a
strengthened form in \cite{Kr4}. The constructions applied in the
proofs essentially used the univalence in the canonical disk
$\D^*$. Theorem 2.1 allows us to solve a similar question 
for the generalized \Gr \ norm $\vk_{D^*}$ of the functions 
univalent in an arbitrary quasidisk $D^*$.

\begin{thm} For any function $f \in \Sigma^0(D^*)$, there exists 
a sequence
$\{f_n\} \subset \Sigma^0(D^*)$ with $\vk_{D^*}(f_n) < k(f_n)$
convergent to $f$ locally uniformly in $D^*$.
\end{thm}

\subsection{}
There is a related conjecture posed in \cite{KK1} that $f \in
\Sigma^0$ with $\vk(f) = k(f)$  cannot be the limit functions
 of locally uniformly convergent sequences $
\{f_n\} \subset \Sigma^0$ with $\vk(f_n) = k(f_n)$. 

Its proof is given in \cite{Kr6}. The main arguments involve a special 
\hol \ motion of the disk and can be appropriately extended to generic 
quasidisks, i.e., to the generalized \Gr \ norm (cf. Section 7).

\subsection{Connection with Fredholm eigenvalues}
The Fredholm eigenvalues $\rho_n$ of a smooth closed Jordan curve
$L \subset\hC$ are the eigenvalues of its double-layer potential,
i.e., of the integral equation
$$
u(z) +  \fc{\rho}{\pi} \int\limits_L \ u(\z)
\fc{\partial}{\partial n_\z} \log \frac{1}{|\z - z|} ds_\z = h(z),
$$
which has many applications. The least positive eigenvalue
$\rho_1 = \rho_L$ plays a crucial role, since by the
K\"{u}hnau-Schiffer theorem it is reciprocal to the \Gr \ norm of
the Riemann mapping function of the exterior domain of $L$. This
value is defined for any oriented closed Jordan curve $L \subset
\hC$ by
$$
\fc{1}{\rho_L} = \sup \ \fc{|\mathcal D_G (u) - \mathcal D_{G^*}
(u)|} {\mathcal D_G (u) + \mathcal D_{G^*} (u)},
$$
where $G$ and $G^*$ are, respectively, the interior and exterior
of $L; \ \mathcal D$ denotes the Dirichlet integral, and the
supremum is taken over all functions $u$ continuous on $\hC$ and
harmonic on $G \cup G^*$.

Until now, no general algorithms exist for finding these values
for the given \qc \ curves. The problem was solved only for some
specific classes of curves, so in general one can use only a rough
estimate for $\rho_L$ by Ahlfors' inequality 
\be\label{2.11}
\fc{1}{\rho_L} \le q_L,
\end{equation}
where $q_L$ is the minimal dilatation of \qc \ reflections across
the given curve $L$, (that is, of the orientation reversing \qc \
homeomorphisms of $\hC$ preserving $L$ point-wise); see, e.g.,  
\cite{Ah2}, \cite{Kr6}, \cite{Ku3}.

Corollary 2.3 provides the following improvement of Ahlfors'
inequality.

\begin{thm} For any quasicirle $L = f(S^1), \ f \in \Sigma^0$,
\be\label{2.12} 
\begin{aligned}\fc{1}{\rho_L} &= \sup \{|\langle
\nu_{S_f}, \psi \rangle_\D|: \ \psi \in A_1^2, \ \|\psi
\|_{A_1(\D) = 1}\} + O(\|S_f\|_\B^2)  \\
&= \sup \{|\langle \mu_0(1/\ov{z})\ov{z}/z, \psi\rangle_\D|: \
\psi \in A_1^2, \ \|\psi \|_{A_1(\D) = 1}\} + O(\|\mu_0\|_\iy^2),
\end{aligned}
\end{equation}
where $\mu_0(z) = g_z/g_{\ov z}$ is the complex dilatation of
extremal quasireflection over the curve $L$ and both remainders are
estimated uniformly for $\|\mu_0\| \le k_0 < 1$.
\end{thm}

\section{Proof of Theorem 2.5}

\noindent $(i)$ \ First observe that the generalized \Gr \
coefficients $\beta_{m n}(f^\mu)$ of the functions $f^\mu \in
\Sigma(D^*)$ generate for each $\x = (x_n) \in l^2$ with $\|\x\| =
1$ the \hol \ maps 
\be\label{3.1} 
h_{\x}(\mu) =
\sum\limits_{m,n=1}^\iy \ \beta_{m n} (f^\mu) x_m x_n : \
\Belt(D)_1 \to \D,
\end{equation}
and $\sup_{\x} |h_{\x}(f^\mu)| = \vk_{D^*}(f^\mu)$.

The holomorphy of these functions follows from the holomorphy of
coefficients $\beta_{m n}$ with respect to \Be \ coefficients $\mu
\in \Belt(D)_1$ mentioned above using the estimate
\be\label{3.2}
\Big\vert \sum\limits_{m=j}^M \sum\limits_{n=l}^N \ \beta_{mn} x_m
x_n \Big\vert^2 \le \sum\limits_{m=j}^M |x_m|^2
\sum\limits_{n=l}^N |x_n|^2
\end{equation} 
which holds for any finite $M, N$ and 
$1 \le j \le M, \ 1 \le l \le N$.
This estimate is a simple corollary of the Milin univalence
theorem (cf. [Mi, p. 193], [Po, p. 61]).

Similar arguments imply that the maps (3.1) regarded as functions
of points $\vp^\mu = S_{f^\mu}$ in the \uTs \ $\T$ (with the basepoint
$D^*$) are \hol \ on $\T$.

Now, let a sequence $\{f_p\} \subset \Sigma^0(D^*)$ be convergent
to $f_0$ uniformly on compact subsets of $\D^*$. Denote their
generalized \Gr \ coefficients by $\beta_{mn}^{(p)}$. Then, for
any $M, N < \iy$ and any fixed $\x = (x_n) \in S(l^2)$,
$$
\Big\vert \sum\limits_1^M \sum\limits_1^N \ \beta_{mn}^{(0)} \ x_m
x_n \Big\vert = \lim\limits_{p\to \iy} \ \Big\vert \sum\limits_1^M
\sum\limits_1^N \ \beta_{mn}^{(p)} \ x_m x_n \Big\vert \le
\liminf\limits_{p\to \iy} \vk_{D^*} (f_p).
$$
Taking the supremum over $\x$ in the left-hand side yields the
desired inequality
$$
\vk_{D^*} (f_0) \le \liminf\limits_{p\to \iy} \vk_{D^*} (f_p).
$$
\bk\noindent $(ii)$ \ Since for any $\mu \in \Belt(D)_1$,
$$
\vk_{D^*}(\vp^\mu) = \sup_{\x\in S(l^2)} |h_{\x}(\vp^\mu)|, \quad
\vp^\mu = S_{f^\mu},
$$
the function $\vk(\vp)$ possesses, together with $h_{\x}(\vp)$,
the mean value inequality property. To get the plurisubharmonicity
of $\vk(\vp)$, one needs to establish its upper semicontinuity.
Using the holomorphy of functions (3.1), one can derive much more.

For any fixed $\x \in S(l^2)$, the function $h_\x(\vp) -
h_\x(\vp_0)$ is a \hol \ map of the ball
$$
\{\vp \in \T: \ \|\vp - \vp_0\|_\B < d\}, \quad d = \dist (\vp_0,
\partial \T)
$$
into the disk $\{|w| < 2\}$. Hence, by Schwarz's lemma,
$$
|h_\x(\vp) - h_\x(\vp_0)| \le \fc{2}{d} \|\vp - \vp_0\|,
$$
and
$$
||h_\x(\vp)| - |h_\x(\vp_0)|| \le |h_\x(\vp) - h_\x(\vp_0)| \le
\fc{2}{d} \|\vp - \vp_0\|.
$$

Now assume that $\vk_{D^*}(\vp) \ge \vk_{D^*}(\vp_0)$ and pick a
maximizing sequence $h_{\x_m}(\vp)$ so that
$$
\lim\limits_{m\to\iy} |h_{\x_m}(\vp)| = \vk_{D^*}(\vp).
$$
Then, since the estimate holds for any $\x \in S(l^2)$, one gets
$$
0 < \vk_{D^*}(\vp) - \vk_{D^*}(\vp_0) \le \vk_{D^*}(\vp)
-\limsup\limits_{m\to \iy} |h_{\x_m}(\vp_0)| \le \fc{2}{d} \|\vp -
\vp_0\|.
$$
In the same way, if $\vk_{D^*}(\vp_0) \ge \vk_{D^*}(\vp)$,
$$
0 < \vk_{D^*}(\vp_0) - \vk_{D^*}(\vp) \le \fc{2}{d} \|\vp -
\vp_0\|,
$$
which implies the Lipschitz continuity of $\vk_{D^*}$ in a
neighborhood of $\vp_0$, completing the proof of the theorem.

\section{Proofs of Theorem 2.1} 

Note that if $\vk_{D^*}(f^\mu) = k(f^\mu) = \|\mu\|_\iy$, 
then 
$$ 
\vk_{D^*}(f^{t\mu}) = k(f^{t\mu}) \quad \text{for all} \ \ |t| < 1. 
$$ 
This follows, for example, from subharmonicity of the function  
$\vk_{D^*}(f^{t\mu})$ in $t$ on the unit disk giving subharmonicity 
of the ratio 
$$ 
g(t) = \fc{\vk_{D^*}(f^{t\mu})}{k(f^{t\mu})} = 
\fc{\vk_{D^*}(f^{t\mu})}{|t|} \quad \text{for all} \ \ |t| < 1. 
$$

We first consider the
case $D^* = \D^*$ which sheds light to key features. For $f \in
\Sigma^0$, the functions (3.1) are of the form 
\be\label{4.1}
h_{\x}(\mu) = \sum\limits_{m,n=1}^\iy \ \sqrt{m n} \ \a_{m
n}(f^\mu) \ x_m x_n : \ \Belt(\D)_1 \to \D.
\end{equation}
Take, for a given function $f$, an extremal coefficient $\mu$ (i.e., 
such that $k(f) = \|\mu\|_\iy$) and consider its extremal disk
$$
\D(\mu) = \{t \mu/\|\mu\|_\iy: \ |t| < 1\} \subset \Belt(\D)_1. 
$$ 
Put $\mu^* = \mu/\|\mu\|_\iy$. 
We apply to $h_{\x}(f^{t\mu^*})$ the well-known improvement of the
classical Schwarz lemma (see \cite{BM}, \cite{Go}) which asserts
that a \hol \ function
$$
g(t) = c_m t^m + c_{m+1} t^{m+1} + \dots: \D \to \D \quad (c_m
\neq 0, \ \ m \ge 1),
$$
in $\D$ is estimated by \be\label{4.2} |g(t)| \le |t|^m \fc{|t| +
|c_m|}{1 + |c_m| |t|},
\end{equation}
and the equality occurs only for
$$
g_0(t) = t^m (t+ c_m)/(1 + \ov c_m t).
$$

To calculate the corresponding constant $\a(f)$ in (2.3), one can
use the variational formula for  $f^\mu(z) = z + b_0 + b_1 z^{-1}
+ \dots \in \Sigma^0$ with extensions satisfying $f^\mu(0) = 0$.
Namely, for small $\|\mu\|_\iy$, 
\be\label{4.3} f^\mu(z) = z -
\fc{1}{\pi} \iint\limits_\D \mu(w) \left( \fc{1}{w - z} -
\fc{1}{w} \right) du dv + O(\|\mu^2\|_\iy), \quad w = u + iv,
\end{equation}
where the ratio $O(\|\mu^2\|_\iy^2)/\|\mu^2\|_\iy^2$ is uniformly
bounded on compact sets of $\C$. Then
$$
b_n = \fc{1}{\pi} \iint\limits_\D \mu(w) w^{n-1} du dv +
O(\|\mu^2\|_\iy), \quad n = 1, 2, \dots,
$$
and from (1.2), \be\label{4.4} \a_{mn}(\mu) = - \pi^{-1}
\iint\limits_\D \mu(z) z^{m+n-2} dx dy + O(\|\mu\|_\iy^2), \quad
\|\mu\|_\iy \to 0.
\end{equation}
Hence, the differential at zero of the corresponding map
$h_\x(t\mu^*)$ with $\x = (x_n) \in S(l^2)$ is given by
\be\label{4.5} 
d h_{x}(0) \mu^* = - \fc{1}{\pi} \iint\limits_\D
\mu^*(z) \sum\limits_{m+n=2}^\iy \sqrt{m n} \ x_m x_n z^{m+n-2} dx
dy.
\end{equation}
On the other hand, as was established in \cite{Kr2}, the elements
of $A_1(\D)^2$ are represented in the form
$$
\psi(z) = \om(z)^2 = \fc{1}{\pi} \sum\limits_{m+n=2}^\iy \
\sqrt{mn} \ x_m x_n z^{m+n-2},
$$
with $\|\x\|_{l^2} = \|\om\|_{L_2}$. Thus, by (4.2), for any $\mu
= t \mu^*$,
$$
|h_\x(\mu)| \le |t| \fc{|t| + |\langle \mu^*, \psi\rangle_\D|}{1 +
|\langle \mu^*, \psi\rangle_\D| |t|},
$$
and $k(f^\mu) = |t|$. Taking the supremum over $\x \in S(l^2)$,
one derives the estimate (2.3).

To analyze the case of equality, observe that if $\a(f^\mu) = 1$,
the second factor in the right-hand side of (2.3) equals $1$, and
this inequality is reduced to $\vk(f^\mu) \le |t| = k(f^\mu)$. But
it was shown in \cite{Kr2} that the equality $\a(f^\mu) = 1$ is
the necessary and sufficient condition to have $\vk(f) = k(f)$.
This completes the proof of the theorem for the canonical disk
$\D^*$.

\bk The case of a generic quasidisk $D$ is investigated along the
same lines using the results established by Milin \cite{Mi} for
the kernels and orthonormal systems in multiply connected domains.
We apply these results to simply connected \qc \ domains $D^*$.
Similar to (4.3),
\be\label{4.6} 
f^\mu(z) = z - \fc{1}{\pi} \iint\limits_D \mu(w) \left( \fc{1}{w
- z} - \fc{1}{w} \right) du dv + O(\|\mu^2\|_\iy),
\end{equation}
but now the kernel of this variational formula is represented 
for $z$ running over a subdomain of $D^*$ bounded by the level
line $G(z, \z) = \rho(w)$ of the Green function of $D^*$ 
in the form
\be\label{4.7} 
\fc{1}{w - z} = \sum\limits_1^\iy P_n^\prime(w) \vp_n(z),
\end{equation}
where $\vp_n = \chi^n$ are given in (1.6) and $P_n$ are
well-defined polynomials; the degree of $P_n$ equals $n$. These
polynomials satisfy 
\be\label{4.8} 
\fc{1}{\pi} \ \iint\limits_D
P_m^\prime(z) P_n^\prime(z) dx dy + \fc{1}{\pi} \
\iint\limits_{D^*} r_m^\prime(z) r_n^\prime(z) dz dy = \dl_{m n},
\end{equation}
where the functions $r_n$ are generated by 
$$
R_{D^*}(z, \z) = \sum\limits_1^\iy r_n(z) \vp_n(\z)
$$ 
(see (1.4))  
and in our case, due to what was mentioned in Section 1.2, 
vanish identically
on $D^*$. Hence, (4.8) assumes the form
$$
\langle P_m^\prime, P_n^\prime \rangle_D = \pi \dl_{m n},
$$
which means that the polynomials $P_n^\prime(z)/\sqrt{\pi}$ 
form an orthonormal system in $A_1^2(D)$. It is
proved in \cite{Mi} that this system is complete.

Noting that for any fixed $z$ the equality (4.7) is 
extended holomorphically to all $w \in D$, one derives 
from (4.6) and (4.8) the following generalization of (4.4).  
From (1.5),  
$$ 
\fc{f^\mu(z) - f^\mu(\z)}{z - \z} 
= 1 - \fc{1}{\pi} \iint\limits_D 
\fc{\mu(w) du dv}{(w - z)(w - \z)} + O(\|\mu^2\|_\iy)
$$  
and 
$$
\begin{aligned}  
- \log \fc{f^\mu(z) - f^\mu(\z)}{z - \z} 
&= - \log \Bigl[ 1 - \fc{1}{\pi} \iint\limits_D 
\fc{\mu(w) du dv}{(w - z)(w - \z)}\Bigr] + O(\|\mu^2\|_\iy) \\ 
&= \fc{1}{\pi} \iint\limits_D 
\fc{\mu(w) du dv}{(w - z)(w - \z)} + O(\|\mu^2\|_\iy) \\  
&= \fc{1}{\pi} \iint\limits_D \mu(w) 
\sum\limits_1^\iy P_m^\prime(w) \vp_m(z) 
\sum\limits_1^\iy P_n^\prime(w) \vp(\z) du dv + O(\|\mu^2\|_\iy),  
\end{aligned}
$$ 
where the ratio $O(\|\mu^2\|_\iy)/\|\mu^2\|_\iy$ is uniformly 
bounded on compact sets of $\C$. 
Comparison with the representation  
$$ 
- \log \fc{f^\mu(z) - f^\mu(\z)}{z - \z} 
= \sum\limits_1^\iy \beta_{m n} \vp_m(z) \vp_n(\z) 
\quad (\vp_n = \chi^n)  
$$ 
yields 
\be\label{4.9}
\wh \beta_{m n}(\mu) = - \fc{1}{\pi} \iint\limits_D \mu(z)
P_m^\prime(w) P_n^\prime(w) du dv + O(\|\mu^2\|_\iy),  
\end{equation} 
which provides the representation of differentials of 
\hol \ functions $\mu \mapsto \wh \beta_{m n}(\mu)$ on 
$\Belt(D)_1$ at the origin. Using the estimate (3.2) 
ensuring the holomorphy of the corresponding functions (3.1) 
on this ball, we get instead of (4.5) that the differential 
of $h_\x(\mu)$ at zero is represented in the form  
\be\label{4.10}
d h_\x(\mathbf 0) \mu^* := \beta_{m n}(f^\mu) = 
- \fc{1}{\pi} \iint\limits_\D \mu^*(z)
\sum\limits_{m,n=1}^\iy  x_m x_n \ P_m^\prime(z) P_n^\prime(z) dx
dy, \quad \x = (x_n) \in S(l^2).
\end{equation} 
Now one can apply the same arguments as in the concluding part of 
the proof in the previous special case and get straightforwardly 
the estimate (2.3) for the general case. 

\bk\noindent 
\textbf{Remark}. \ 
The equality (4.5) yields that in the case $D^* = \D^*$ 
the constant (2.2) for every $f \in \Sigma^0(D^*)$ is represented 
in the form 
\be\label{4.11}
\a_D(f) = \sup_{\x=(x_n)\in S(l^2)} \ 
\fc{1}{\pi \|\mu\|_\iy} \Big\vert \iint\limits_{|z|<1} 
\mu(z) \ \sum\limits_{m+n\ge 2}^\iy  
\sqrt{m n} \ x_m x_n z^{m+n-2} dxdy \Big\vert,     
\end{equation} 
where $\mu$ is any extremal \Be \ coefficient in the 
equivalence class $[f]$.

\section{Proof of Theorem 2.2} 
 
It follows from the proof of Theorem 2.1 that 
$\vk_{D^*}(f) = k(f)$ if and only if $\a_D(f) = 1$. 
So it remains to establish the equlity (2.5), 
provided that the extremal extension of $f$ to $D$ is 
of \Te \ type, with \Be \ coefficient 
$\mu_0 = k |\psi_0|/\psi_0, \ k = k(f)$.  

Pick a sequence $\{\psi_p = \om_p^2\} \subset A_1^2(D)$ with
$\|\om_p\|_{L_2(D)} = 1$ for which 
$$
\lim\limits_{p\to \iy} |\langle |\psi_0|/\psi_0, 
\psi_p \rangle_D| = 1. 
$$ 
This sequence is convergent uniformly on compact sets in $D$ to 
a \hol \ function $\vp \in A_1^2$.

If $\vp(z) \equiv 0$, the sequence $\{\psi_p\}$ should be
degenerate for the coefficient $|\psi_0|/\psi_0$, which is
impossible for \Te \ extremal coefficients. Thus 
$\vp \ne \mathbf 0$, and 
\be\label{5.1} 
|\langle |\psi_0|/\psi_0, \vp \rangle_D|
\le \lim\limits_{p\to \iy} |\langle |\psi_0|/\psi_0, \psi_p
\rangle_D| = 1.
\end{equation}
It remains to show that, under assumptions of the theorem, the
left inequality in (5.1) must be an equality (hence 
$\psi_0 = \vp$).  
We may assume that $f(0) = 0$ (passing if needed to 
$f_1(z) = f(z) - f(0)$).

Noting that in view of (2.2) and (4.10) each $\wt \psi_p$ 
is represented in the form
$$
\psi_p(z) = \fc{1}{\pi} \sum\limits_{m,n=1}^\iy \ 
x_m^{(p)} x_n^{(p)} P_m^\prime(z) P_n^\prime(z)  
\quad \text{with} \ \ \x^{(p)} = (x_n^{(p)}) \in S(l^2) 
$$ 
and selecting if needed a subsequence from $\x^{(p)}$ convergent in
$l^2$ to $\x^{(0)} = (\x_n^{(0)})$, one gets  
$\lim\limits_{p\to \iy} \x_n^{(p)} = \x_n^{(0)}$ for each $n \ge 1$, 
and by the above remark $\x^{(0)} \ne \mathbf 0$. 
This implies that $\vp$ as the weak limit of $\psi_p$ is of the 
form  
$$ 
\vp(z) = \pi^{-1} \sum\limits_{m,n=1}^\iy \ \ x_m^{(0)} x_n^{(0)}
P_m^\prime(z) P_n^\prime(z),  
$$
and $\|\x^0\|_{l^2} = 1$ (in view of maximality of
$\vk_{D^*}(f^{k|\vp|/vp})$.
The variation (4.9) yields that the \Gr \ coefficients 
of $f^{t\mu_0} = f$ and $f^{t|\vp|/\vp}$ are related by  
$$
\beta_{m n}(f^{t \mu_0}) = \beta_{m n}(f^{t|\vp|/\vp}) + O(t^2), 
\quad t \to 0, 
$$
and, letting $t \to 0$, 
$$
\langle |\psi_0|/\psi_0, P_m^\prime P_n^\prime\rangle_\D 
= \langle |\vp|/\vp, P_m^\prime P_n^\prime\rangle_\D 
\quad \text{for all} \ \ m,n \ge 1. 
$$
Extension of these functionals to $A_1(D)$ by Hahn-Banach yields 
$$
\langle |\psi_0|/\psi_0 - |\vp_1|/\vp_1, \psi \rangle_\D =
0 \quad \text{for any} \ \ \psi \in A_1(\D).
$$
As is well known (see, e.g., \cite{GL}, \cite{Kr1}), such equality 
is impossible for the \Te \ extremal coefficients unless 
$\psi_0 = \wt \vp$. This completes the proof of the theorem.

\section{Proof of Corollary 2.4 and of Theorem 2.7}

\bk\noindent
\textbf{Proof of Corollary 2.4}. Theorem 2.1 and the inequality (2.6) 
(for $D^* = \D^*$) yield that the equality (2.1) holds for all 
$f \in \Sigma^0$
admitting the \Te \ extremal extensions $f^{\mu_0}$ across $S^1$.
The Schwarzians derivatives $S_{f^{\mu_0}}$ of such $f$ are 
Strebel's points of the space $\T$ (cf., e.g.
\cite{GL}, \cite{St}).

For sufficiently small $|t|$, the Schwarzians $\vp_t =
S_{f^{t\mu_0^*}}$ determine by (2.8) the harmonic \Be \
coefficients of the Ahlfors-Weill extension of the maps 
$f^{t\mu_0^*}$ across the unit circle $S^1 = \partial \D^*$.  
In view of the
characteristic property of extremal \Be \ differentials, we have
for any such $\mu_0^*$ the equality
$$
\nu_{\vp_0} = t\mu_0^* + \sigma_0, \quad  \sigma \in A_1(\D)^\bot,
$$
where
$$
A_1(\D)^\bot = \{\nu \in \Belt(\D)_1: \ \langle \nu,\psi\rangle_\D
= 0 \ \text{for all} \ \psi \in A_1(\D)\}
$$
is the set of infinitesimally trivial \Be \ coefficients
(see e.g. \cite{GL}, \cite{Kr1}).

Since, due to \cite{GL}, the set of Strebel's points are open 
and dense in \Te \ spaces,
the equality (2.9) (and its equivalent(2.7)) must hold for all
points $\vp = S_f$ (with sufficiently small norms), which
completes the proof of the corollary.

Note that by the same reasons the inequality (2.4) holds for
all $f \in \Sigma^0$.

\bk
\noindent\textbf{Proof of Theorem 2.7}. 
Since all quantities in (2.11) are invariant under the action of
the M\"{o}bius group $PSL(2, \C)/\pm \mathbf{1}$, it suffices to
use \qc \ homeomorphisms $f$ of the sphere $\hC$ carrying the unit  
circle $S^1$ onto $L$ whose \Be \ coefficients
$\mu_f(z) = \partial_{\ov{z}} f/\partial_z f$ are supported in the
unit disk $\D$ and which are hydrodynamically normalized near the
infinite point, i.e., with restrictions $f|\D^* \in \Sigma^0$. 
Then the reflection coefficient $q_L$ equals the minimal
dilatation $k(w^\mu) = \|\mu\|_\iy$ of \qc \ extensions $w^\mu$ of
$f|\D^*$ to $\hC$, and Theorem 2.7 immediately follows from
Corollaries 2.3 and 2.4.

\section{Proof of Theorem 2.6}

Again, in view of density of Strebel's points in $\T$, 
it suffices to prove this theorem for $f \in \Sigma^0(D)$ with 
\Te \ extensions $f^{k|\psi|/\psi}$ to $\hC$ defined by 
quadratic differentials $\psi \in A_1(D)$ of the form 
$$ 
\psi(z) = c_m z^m + O(z^{m+1}) \quad \text{near} \ \ z = 0 \ \ 
(m \ge 0).   
$$  
Take
$$
\psi_n(z) = \psi(z) + \fc{c_1^{(n)}}{z}
$$
with $c_1^{(n)} \to 0$ as $n \to \iy$, and consider the maps
$f_n(z) = f^{k|\psi_n|/\psi_n}$. These maps are convergent to 
$f(z)$ uniformly on compact sets of $\C$ and $k(f_n) = k$.  
However, since every $\psi_n$ has a simple pole at the origin, 
$$
\a_{D^*}(f_n) = \sup \ \Big\{ \Big\vert \iint\limits_D
\fc{|\psi_n(z)|}{\psi_n(z)} \vp(z) dx dy \Big\vert: \ \vp \in
A_1^2(D), \ \|\vp\|_{A_1} = 1 \Big\} < 1.
$$
Hence, by (2.3), for any $n$,
$$
\vk_{D^*}(f_n) \le k \fc{k + \a(f)}{1 + \a(f) k} < k, 
$$
completing the proof.

\section{Examples} 

\subsection{} 
It follows from Theorem 2.2 (equality (2.5)) that 
$\vk_{D^*}(f) < k(f)$
{\em for any quasidisk $D^* \ni \iy$ and any $f \in \Sigma^0(D^*)$ 
having the \Te \ extension $f^\mu$ to $D$ with    
$\mu = k |\psi_0|/\psi_0$, where $\psi_0$ is \hol \ and 
has zeros of odd order in} $D$. 
The simplest example of such $\psi_0$ is given by 
$\psi_0(z) = z^p$ with an odd integer $p \ge 1$. 

To get other examples, one can pick 
$\psi_0 = g(z)^p$, where $g(z)$ is a conformal map of $D$ 
onto the unit disk with $g(0) = 0, \ g^\prime(0) > 0$.  

An explicit construction of the Riemann mapping 
functions of simply connected domains is a very difficult 
problem. Their representation is known only for some special 
domains. 

\subsection{}
For example, if $D_{\mathcal E}^*$ is the exterior of the ellipse 
$\mathcal E$ with the foci at $-1, 1$ and semiaxes 
$a, b \ (a > b)$, then the branch of the function 
$$ 
\chi(z) = (z + \sqrt{z^2 - 1})/(a + b) 
$$ 
positive for real $z > 1$ maps this exterior onto $\D^*$. 
A conformal map of the interior of this ellipse $D_{\mathcal E}$ 
onto the disk involves an elliptic function. 

As is well known (see \cite{Ne}), an orthonormal basis in the space 
$$ 
A_2(D_{\mathcal E}) = \{\om \in L_2(D_{\mathcal E}): \ 
\om \ \ \text{\hol \ in} \ \ D_{\mathcal E}\} 
$$
is formed by the polynomials 
$$ 
P_n(z) = 2 \sqrt{\fc{n + 1}{\pi}} \ (r^{n+1} - r^{-n-1}) \ U_n(z), 
$$ 
where $r = (a + b)^2$ and $U_n(z)$ are the Chebyshev polynomials 
of the second kind, 
$$ 
U_n(z) = \fc{1}{\sqrt{1 - z^2}} \ \sin [(n + 1) \arccos z], 
\quad n = 0, 1, \dots \ . 
$$ 
Using the Riesz-Fisher theorem, one obtains that each function 
$\psi \in A_2(D_{\mathcal E})$ is of the form (cf. \cite{Kr2}) 
$$ 
\psi(z) = \sum\limits_0^\iy x_n P_n(z), \quad \x = (x_n) \in l^2, 
$$ 
with $\|\psi\|_{A_2} = \|\x\|_{l^2}$.  

By Theorem 2.2, a function $f \in \Sigma^0(D_{\mathcal E}^*)$ 
with \Te \ extension $f^\mu$ to $D_{\mathcal E}$ satisfies 
\be\label{8.1} 
\vk_{D_{\mathcal E}^*}(f) = k(f) 
\end{equation}
if and only if 
$$ 
\mu(z) = k \ov{\sum\limits_0^\iy x_n^0 P_n(z)} \Big/  
\sum\limits_0^\iy x_n^0 P_n(z) 
$$ 
with some $\x^0 = (x_n^0) \in S(l^2)$. 
More generally, a function $f \in \Sigma^0(D_{\mathcal E}^*)$ 
obeys (8.1) if and only if any its extremal \Be \ coefficient 
$\mu \in \Belt(D_{\mathcal E})$ satisfies 
$$ 
\sup \Big\vert\Big\langle 
\mu, \sum\limits_{m,n\ge 0}^\iy  
x_m x_n P_m P_n \Big\rangle_{D_{\mathcal E}} \Big\vert 
= \|\mu\|_\iy,   
$$ 
taking the supremum over all $\x = (x_n) \in l^2$ with $\|\x\| = 1$. 
Note also that for every $f \in \Sigma^0(D_{\mathcal E}^*)$, its 
constant $\a_{D_{\mathcal E}}(f)$ is given explicitly by 
$$ 
\a_{D_{\mathcal E}}(f) = \sup_{\x=(x_n)\in S(l^2)} \ 
\Big\vert \iint\limits_{D_{\mathcal E}} 
\fc{\mu(z)}{\|\mu\|_\iy} \ \sum\limits_{m,n\ge 0}^\iy  
x_m x_n P_m(z) P_n(z) dx dy \Big\vert,     
$$ 
taking any extremal $\mu$ in the equivalence class $[f]$.   

\subsection {} 
The expansion (1.6) contains a conformal map $\chi: \ D^* \to \D^*$, 
while the basic quantity $\a_D(f)$ is connected with conformal 
maps of the complementary quasidisk $D$. 
The only known non-trivial example with a simple connection 
between these maps is the Cassini curve 
$L = \{z: \ |z^2 - 1| = c\}$ with $c > 1$. It is given in \cite{HK}.  
Here $\chi^{-1}(z) = \sqrt{1 + c z^2}$,   
and the branch of 
\be\label{8.2} 
g(z) = z \sqrt{(c^2 -1)/(c - z^2)} 
\end{equation} 
maps conformally the unit disk onto the interior of $L$ 
with $g(0) = 0, \ g^\prime(0) > 0$.  

Using the function (8.2), one gets that for every univalent 
function $f(z)$ in the domain $D^* = \{|z^2 - 1| > c\}$ 
with hydrodynamical normalization, its constant $\a_D(f)$ is 
given, due to (4.11), by 
$$ 
\a_D(f) = \sup_{\x\in S(l^2)} \ 
\fc{1}{\pi \|\mu\|_\iy} \Big\vert \iint\limits_{|z|<1} 
\mu \circ g(z) \fc{\ov{g^\prime(z)}}{g^\prime(z)} \ 
\sum\limits_{m+n\ge 2}^\iy  
\sqrt{m n} \ x_m x_n z^{m+n-2} dxdy \Big\vert,     
$$ 
taking again an extremal \Be \ coefficient $\mu$ in the class 
$[f]$.

\section{Grunsky norm and complex homotopy}  

Every function $f \in \Sigma(D^*)$ generates a \hol \ homotopy by 
\be\label{9.1} 
f(z, t) = f_t(z) := t f \circ g^{-1}[t g(z)]: 
\ D^* \times \D \to \hC, 
\end{equation} 
where $g$ maps conformally $D^*$ onto $D^*$ with 
$g(\iy) = \iy, \ g^\prime(\iy) > 0$. This homotopy satisfies 
$f(z, 0) = z, \ f(z, 1) = f(z)$
and $ f(z, t) = z + \wt b_0 t + \wt b_1 t^2 z^{-1} + \dots$ 
near $z = \iy$. 
The curves $\{z = - \log |t|\}$ are the level lines of Green's 
function $g_{D^*}(z, \iy) = - \log |g(z)|$ of $D^*$. 

Consider the Schwarzians $S_{f_t}(z) = S_f(z, t)$. Then the map 
$t \mapsto S_f(\cdot, t)$ is \hol \ in $t$ for any $z \in D^*$  
and, due to the well-known properties of the functions with 
sup mnorm depending holomorphically on complex parameters, 
this pointwise map induces a \hol \ map 
\be\label{9.2} 
\chi_f: \ t \mapsto S_f(\cdot, t), \quad  
\chi_f(t) = \chi_f(0) + t \chi_f^\prime(0) + \dots, \ \ 
\chi_f(0) = S_{g^{-1}}, 
\end{equation}
of the disk $\{|t| < 1\}$ into the space $\T$.    
We call a level $r = |t| > 0$ \textbf{noncritical} if 
$\chi_f^\prime(r e^{i \theta}) \ne \mathbf 0$ 
for any $\theta \in [0, 2 \pi]$. 
If $\chi_f^\prime(t_0) = \mathbf 0$ then 
$\chi_{f_\eta}^\prime(r) = \mathbf 0$ for 
$\eta = - t_0/|t_0|$.  
In the simplest case of the disk $\D^*$, 
$$  
f_t(z) = t f(z/t) = z + b_0 t + b_1 t^2 z^{-1} + \dots 
$$  
and $S_{f_t}(z) = t^{-2} S_f(z/t)$   
for all $|t| < 1$; then the map (9.1) takes the form   
$$
\chi_f(t) = \fc{\chi_f^{\prime\prime}(0)}{2!} t^2 
+ \fc{\chi_f^{\prime\prime\prime}(0)}{3!} t^3 + \dots \ ,  
$$ 
and the \Gr \ coefficients of $f_t$ are homotopically homogeneous:       
\be\label{9.3}  
\quad \a_{m n} (f_t) = \a_{m n} (f) \ t^{m+n} 
\quad \text{for all} \quad m, n \ge 1.   
\end{equation}  
The homotopy disk 
$$ 
\D(S_f) := \chi_f(\D) = \{S_{f_t}: \ |t| < 1\}
$$ 
has cuspidal singularities in the critical 
points of $\chi_f$. 

For any quasidisk $D^*$ containing the infinite point, we have 

\begin{thm} Let the homotopy function $f_r(z)$ of  
$f\in \Sigma(D^*)$ given by (9.1) satisfy  
\be\label{9.4}
\vk_{D^*}(f_\rho) = k(f_\rho)
\end{equation}
for a noncritical level $\rho \in (0, 1)$. Then  
\be\label{9.5}
\vk_{D^*}(f_r) = k(f_r) \quad \text{for all} \ \ r < \rho. 
\end{equation}
\end{thm}

This theorem answers some questions stated by R. K\"{u}hnau 
in \cite{KK2}. It also has some other interesting applications. 
Apart from some special cases, there is no connection 
between the defining \hol \ quadratic differentials $\psi$ 
and $\psi_r$ of a map $f$ and its homotopies $f_r$. 
Theorems 2.2 and 9.1 give the conditions ensuring the evenness 
of zeroes of $\psi$ and $\psi_r$ (cf. \cite{Kr7}). 

\bk
The proof of Theorem 9.1 essentially involves the curvature 
properties of the \Ko \ metric of \uTs \ $\T$. 
We first recall some background facts underlying the proof.  

We shall use the following strengthening of the fundamental 
Royden-Gardiner theorem given in \cite{Kr3}.  

\begin{prop} The differential (infinitesimal) 
\Ko \ metric $\mathcal K_\T(\vp, v)$
on the tangent bundle $\mathcal T(\T)$ of the \uTs \ $\T$ is
logarithmically \psh \ in $\vp \in \T$, equals the canonical
Finsler structure $F_\T(\vp, v)$ on $\mathcal T(\T)$ generating
the \Te \ metric of $\T$ and has constant \hol \ sectional
curvature $\kappa_{\mathcal K}(\vp, v) = - 4$ on $\mathcal T(\T)$.
\end{prop}

The \textbf{generalized Gaussian curvature} $\kappa_\ld$ of 
an upper semicontinuous Finsler metric $ds = \ld(t) |dt|$ in 
a domain $\Om \subset \C$ is defined by
\be\label{9.6}
\kappa_\ld (t) = - \fc{\mathbf{\D} \log
\ld(t)}{\ld(t)^2},
\end{equation}
where $\mathbf{\D}$ is the \textbf{generalized Laplacian}
$$
\mathbf{\D} \ld(t) = 4 \liminf\limits_{r \to 0} \frac{1}{r^2}
\Big\{ \frac{1}{2 \pi} \int_0^{2\pi} \ld(t + re^{i \theta}) d
\theta - \ld(t) \Big\}
$$
(provided that $- \iy \le \ld(t) < \iy$). Similar to $C^2$
functions, for which $\mathbf{\D}$ coincides with the usual
Laplacian, one obtains that $\ld$ is \sh \ on $\Om$ if and only if
$\mathbf{\D} \ld(t) \ge 0$; hence, at the points $t_0$ of local
maximuma of $\ld$ with $\ld(t_0) > - \iy$, we have $\mathbf{\D}
\ld(t_0) \le 0$.

The sectional \textbf{\hol \ curvature} of a Finsler metric on a
complex Banach manifold $X$ is defined in a similar way as the
supremum of the curvatures (9.6) over appropriate collections of
\hol \ maps from the disk into $X$ for a given tangent direction
in the image.
The \hol \ curvature of the \Ko \ metric $\mathcal K_X(x, v)$ of any
complete hyperbolic manifold $X$ satisfies $\kp_{\mathcal K} \ge
- 4$ at all points $(x, v)$ of the tangent bundle $\mathcal T(X)$
of $X$, and for the \Ca \ metric $\mathcal C_X$ we have
$\kp_{\mathcal C}(x, v) \le - 4$ (cf., e.g., \cite{AP}, \cite{Di}, 
\cite{Ko}). 

It was istablished in \cite{EE} that the metric 
$\mathcal K_\T(\vp, v) = F_\T(\vp, v)$ is Lipschitz continuous 
on $\T$ (in its Bers' embedding).  

We shall deal with \sh \ circularly symmetric (radial) metrics 
$\ld(t) |dt|$ on a disk $\{|t| < a\}$, i.e., such that 
$\ld(t) = \ld(|t|)$. 
Any such function $\ld(t)$ is monotone increasing in 
$r = |t|$ on $[0, a]$ and convex with respect to $\log r$,  
has one-sided derivatives for each $r < a$ (in particular 
$u^\prime(0) \ge 0$), 
and $r u^\prime(r)$ is monotone increasing (see, e.g., \cite{Ro}).    

\bk\noindent 
\textbf{Proof of Theorem 9.1}. First consider the more simple case 
of the circular disk $D^* = \D^*$ which we use to illustrate 
the main ideas. 

The relations (1.4), (9.3), (9.4) imply that the stretching 
$f_\rho$ 
possesses a \Te \ extension to $\D$ defined by a quadratic 
differential $\psi_\rho \in A_1^2(\D)$ so that  
$\mu_{f_\rho}(z) = k(f_\rho) |\psi_\rho(z)|/\psi_\rho(z)$   
for $|z| < 1$) and 
$$ 
\vk(f_\rho) = k(f_\rho) = \Big\vert \sum\limits_{m,n=1}^\iy \ 
\sqrt{mn} \ \a_{mn}(f) \rho^{m+n} x_m^0 x_n^0 \Big\vert 
$$  
(this common value is attained on some point 
$\x^0 = (x_n^0) \in S(l^2)$). 
Indeed, the corresponding function (3.1) for this $\x^0$ 
(with $\beta_{m n} = \sqrt{m n} \ \a_{m n}$) being restricted 
to the disk $\D(S_f)$ assumes the form  
\be\label{9.7}  
\wt h_{\x^0}(t) = \sum\limits_{m,n=1}^\iy \ 
\sqrt{mn} \ \a_{mn}(f) x_m^0 x_n^0 t^{m+n},      
\end{equation} 
and by(9.4), 
\be\label{9.8} 
|\wt h_{\x^0}(\rho)| = \vk(f_\rho) = k(f_\rho). 
\end{equation}  
The series (9.7) defines a \hol \ selfmap of the disk 
$\{|t| < 1\}$.   

Noting that the homotopy $f(z, t)$ is a \hol \ motion of 
the disk $\D^*$ parametrized by $t \in \D$ and applying 
to it the basic lambda-lemma for these motions, one obtains 
that each fiber map $f_t(z) = f(z, t)$ extends to a \qc \ 
automorphism of the whole sphere $\hC$ so that the 
\Be \ coefficient 
$\mu(z, t) = \ov{\partial} f_t/\partial f_t \in \Belt(\D)_1$ 
is a $L_\iy$-\hol \ function of $t \in \D$ (and generically 
not extremal).  
If the derivative of the map (9.2) vanishes at some point $t_0$, 
$\chi_f^\prime(t_0) = \mathbf 0$, then also 
$\fc{d}{dt} \mu(z, t)|_{t=t_0} = \mathbf 0$, and the \hol \ 
dependence of the function (9.7) on $S_{f_t}$ 
and on $\mu(\cdot, t)$ implies    
$$ 
\wt h_{\x^0}^\prime(t_0) = 0.
$$ 
Hence, all critical points of the map (9.2) are simultaneously 
critical for the function (9.7) (though $\wt h_{\x^0}$ 
can have extra critical points which are regular for $\chi_f$).  

We apply the functions (3.1) to the explicit construction of   
some \sh \  
Finsler metrics on \hol \ disks $\Om = g(\D) \subset \T$, 
pulling back the hyperbolic metric  
$\ld_\D(t)|dt| = |dt|/(1 - |t|^2)$ of $\D$ (assuming that the 
\Gr \ coefficients $\a_{m n}$) are given).  
In fact, we shall use these metrics only on the homotopy 
disk $\D(S_f)$ and on geodesic \Te \ disks passing through 
the origin and points of $\D(S_f)$. These metrics 
are dominated by the \Ko-\Te \ metric of the space $\T$. 
The functions  
$$ 
h_{\x,g} (t) := h_\x(S_f \circ g(t)) 
= \sum\limits_{m,n=1}^\iy \ 
\sqrt{mn} \ \a_{mn}(S_f \circ g(t)) x_m x_n,   
\quad \x \in S(l^2), 
$$ 
define  \hol \ maps $\D \to \Om \to \D$ and conformal metrics  
$\ld_{h_{\x,g}}(t) |dt|$ with   
$$
\ld_{h_{\x,g}}(t) 
= |h_{\x,g}^\prime (t)|/(1 - |h_{\x,g} (t)|^2), \quad t \in \D       
$$
of Gaussian curvature $-4$ at noncrical points. 
We take the upper envelope of these metrics 
\be\label{9.9}
\wt \ld_\vk(t) = \sup \{\ld_{h_{\x,g}}(t): \ \x \in S(l^2)\} 
\end{equation}
and its upper semicontinuous regularization
$$
\ld_\vk(t) =\limsup\limits_{t'\to t} \wt \ld_\vk(t'), 
$$ 
getting a logarithmically \sh \ metric on $\Om$.   
In fact, one can show, similarly to Theorem 2.5, 
that this regularization does not 
change (increase) $\ld_\vk$, i.e. $\ld_\vk = \wt \ld_\vk$, 
in view of continuity. 

\bk
Now recall that a conformal metric $\ld_0(t) |dt|$ 
is called \textbf{supporting} for $\ld(t)|dt|$ 
at a point $t_0$ if $\ld_\vk(t_0) = \ld_0(t_0)$ and 
$\ld_0(t) < \ld_\vk(t)$ for all $t \setminus \{t_0\}$ 
from a neighborhood of $t_0$.     

\begin{lem} If a conformal metric $\ld$ in a domain 
$\Om$ has at any its noncritical point $t_0$ a 
supporting \sh \ metric $\ld_0$ of Gaussian curvature 
at most $- 4$, then $\ld$  is \sh \ on $\Om$ and 
its generalized Gaussian curvature also is at most $- 4$  
in all noncritical points. 
\end{lem}   

\medskip 
\noindent 
\textbf{Proof}. Since the space $\B(D^*)$ is dual to $A_1(D^*)$, 
the sequences $\{h_\x(\vp)\}$ are  convergent, by the 
Alaoglu-Bourbaki theorem, in weak* topology to \hol \ 
functions $\T \to \D$. 
This yields that the metric (9.9) has a supporting 
metric $\ld_0(t)$  
in a neighborhood $U_0$ of any noncritical point $t_0 \in \D$, 
which means that $\ld_\vk(t_0) = \ld_0(t_0)$ and 
$\ld_0(t) < \ld_\vk(t)$ for all $t \setminus \{t_0\}$. 
Hence, for sufficiently small $r > 0$,  
$$ 
\frac{1}{r^2} \Bigl(\frac{1}{2 \pi} \int\limits_0^{2\pi} 
\log \ld_\vk(t_0 + r e^{i\theta}) d \theta - \ld_\vk(t_0)\Bigr) 
\ge \frac{1}{r^2} \Bigl(\frac{1}{2 \pi} \int\limits_0^{2\pi} 
\log \ld_0(t_0 + r e^{i\theta}) d \theta - \ld_0(t_0) \Bigr),   
$$
and 
$\mathbf{\D} \log \ld_\vk(t_0) \ge \mathbf{\D} \log \ld_0(t_0)$. 
Since $\ld_\vk(t_0) = \ld_0(t_0)$, one gets  
$$ 
- \frac{\mathbf{\D} \log \ld_\vk(t_0)}{\ld_\vk(t_0)^2} 
\le - \frac{\mathbf{\D} \log \ld_0(t_0)}{\ld_0(t_0)^2} \le - 4,    
$$ 
which completes the proof of the lemma. 

Note that the inequality $\kappa_{\ld} \le - 4$ is equivalent to   
$$
\mathbf{\D} \log \ld \ge 4 \ld^2,
$$ 
where $\mathbf{\D}$ again means the generalized Laplacian. 
Letting $u = \log \ld$, one gets 
$\mathbf{\D} u \ge 4 e^{2 u}$. 

In particular, {\em all this holds for the metrics $\ld_\vk$} 
(cf. \cite{Kr4}).  
Indeed, the space $\B(D^*)$ is dual to $A_1(D^*)$, thus 
by the Alaoglu-Bourbaki theorem  the family 
$\{h_(\vp)\}, \ \vp = S_f \in \T$  
is compact in weak$^*$ topology. The limit functions 
of its subsequences are \hol \ maps of $\T$ and into 
the unit disk.  
This yields that each of the metrics (8.9) has a supporting 
metric in a neighborhood $U_0$ of any noncritical point 
$t_0 \in \D$. 

\bk 
We proceed to the proof of the theorem and note that 
in the case $\Om = \D(S_f)$ the 
enveloping metric (9.9) and both norms $\vk(f_t)$ and $k(f_t)$ 
are circularly symmetric in $t$. 
We determine on this disk also another circularly symmetric 
\sh \ conformal metric majorated by $\ld_\vk$. 

Namely, the map (9.7) generates the metric 
\be\label{9.10} 
\ld_{\wt h_{\x^0}}(t) 
= |\wt h_{\x_0}^\prime (t)|/(1 - |\wt h_{\x^0} (t)|^2)   
\end{equation} 
of Gaussian curvature $- 4$ on $\D$ (again at noncritical points), 
which is supporting for $\ld_\vk(t)$ at $t = \rho$. 
Replacing $\x^0$ by the points 
$\x_\epsilon^0 = (\epsilon x_n^0) \in S(l^2)$ with $|\epsilon| = 1$, 
one gets the corresponding \sh \ metrics 
$\ld_{\wt h_{\epsilon \x^0}}(t) 
= |\wt h_{\epsilon \x_0}^\prime (t)|/(1 - 
|\wt h_{\epsilon \x^0} (t)|^2)$. 
Take their envelope   
\be\label{9.11}  
\ld_0(t) := \sup_\epsilon \ld_{\wt h_{\epsilon \x^0}}(t);  
\end{equation} 
its curvature also is at most $- 4$ in both supporting and 
\hol \ senses. 

\bk 
Our goal now is to prove the equlity 
\be\label{9.12}
\ld_\vk(t) = \ld_d(S_{f_t}, v), 
\end{equation} 
where $\ld_d$ is the restriction to $\D(S_f)$ of the infinitesimal 
\Ko-\Te \ metric on the space $\T$ and $v$ is a tangent vector 
to the \Te \ disk touching $\D(S_f)$ at the point $t$. 
We apply Minda's maximum principle given by  

\begin{lem} \cite{Min} If a function
$u : \ D \to [- \iy, + \iy)$ is upper semicontinuous in a domain
$\Om \subset \C$ and its generalized Laplacian satisfies the
inequality $\mathbf{\D} u(z) \ge K u(z)$ with some positive
constant $K$ at any point $z \in D$, where $u(z) > - \iy$, and
if
$$
\limsup\limits_{z \to \z} u(z) \leq 0 \ \ \text{for all} \
\z \in \partial D,
$$
then either $u(z) < 0$ for all $z \in D$ or else $u(z) = 0$ for all
$z \in \Om$.
\end{lem} 

First observe that  
\be\label{9.13} 
\ld_\vk(\rho) = \ld_d(S_{f_\rho}, v),  
\end{equation} 
which follows from the reconstruction lemma for \Gr \ norm.   

\begin{lem} \cite{Kr4} On any extremal \Te \ disk
$\D(\mu_0) = \{\phi_\T(t \mu_0): \ t \in\D\}$
(and its isometric images in $\T$), we have the equality
\be\label{9.14}
\tanh^{-1}[\vk(f^{r\mu_0})] = \int\limits_0^r \ld_\vk(t) dt.
\end{equation}
\end{lem} 

Indeed, assuming $\ld_\vk(\rho) < \ld_d(S_{f_\rho}, v)$, 
one would have from semicontinuiuty of both sides that such 
strong inequality must hold in a neighborhood of $S_{f_\rho}$ 
in $\T$, but this violates the 
equalities (9.4) and (9.14) for $r = \rho$ (along the corresponding 
\Te \ disk). This proves (9.13).   

The equality (9.13) yields that each of the metrics (9.11), (9.12) 
and $\ld_d(S_{f_t}, v)$ is supported at $t = \rho$ 
by the same metric (9.10). 
Take the annulus 
$\mathcal A_{r_1, r_2} = \{r_1 < \rho < r_2\}$ 
with $r_1 < \rho < r_2$, which does not contain  
the critical points of function (9.7), and put  
$$ 
M = \{\sup \ld_d(t): t \in \mathcal A_{r_1, r_2}\};      
$$ 
then
$\ld_d(t) + \ld_0(t) \le 2M$. Consider on this annulus the function 
$$ 
u(r) = \log \fc{\ld_0(r)}{\ld_d(r)}. 
$$ 
Then (cf. \cite{Min}, \cite{Kr4}), 
$$ 
\D u(r) = \log \ld_0(r) - \ld_d(r) = 4 (\ld_0^2(r) - \ld_d^2(r)) 
\ge 8M (\ld_0(r) - \ld_d(r)), 
$$  
and the elementary estimate 
$M \log(t/s) \ge t - s$ for $0 < s \le t < M$ 
(with equality only for $t = s$) implies  
$$ 
M \log \fc{\ld_0(r)}{\ld_d(r)} \ge \ld_0(r) - \ld_d(r),  
$$ 
and hence, 
$\D u(t) \ge 4 M^2 u(t)$.  

One can apply Lemma 9.4 which implies, in view of the 
equality (9.13), that 
$\ld_0(r) = \ld_\vk(r) = \ld_d(r)$ for all 
$r \in [r_1, r_2]$ (equivalently, $\vk_f(r) = k_f(r)$).  

Now one can fix $\rho < r^\prime < r_2$ and compare 
the metrics $\ld_\vk$ and $\ld_d$ on the disk 
$\{|t| < r^\prime\}$ in a similar way, which yields 
the desired equalities $\ld_\vk(r) = \ld_d(r)$ and 
$\vk_f(r) = k_f(r)$ for all $r \le \rho$, completing 
the proof for the disk $\D^*$.    

The proof for the functions on generic quasidisks $D^*$ 
follows the same lines using the homotopy (9.2).

\bk

\bk
\medskip
{\small\em{ \leftline{Department of Mathematics, Bar-Ilan
University} \leftline{5290002 Ramat-Gan, Israel} \leftline{and
Department of Mathematics, University of Virginia,}
\leftline{Charlottesville, VA 22904-4137, USA}}}

\end{document}